\documentclass{ws-procs9x6}

\begin{document}

\title{A Stochastic Heisenberg Inequality}

\author{C.~Mueller\footnote{\uppercase{S}upported by 
\uppercase{NSA} and \uppercase{NSF} grants.}} 
\address{Department of Mathematics,\\
University of Rochester,\\ 
Rochester, NY 14627, USA\\ 
E-mail: cmlr@math.rochester.edu}

\author{A.~Stan\footnote{\uppercase{S}upported by 
 an \uppercase{NSF} grant.}}

\address{Department of Mathematics,\\
University of Rochester,\\ 
Rochester, NY 14627, USA\\ 
E-mail: astan@math.rochester.edu}  


\maketitle

\abstracts{An analogue of the Fourier transform will be introduced for 
all square integrable continuous martingale processes whose quadratic 
variation is deterministic. Using this transform we will formulate and 
prove a stochastic Heisenberg inequality.}

\section{Introduction}
The Heisenberg inequality in ${\mathbb R}$ (see Folland and 
Sitaram\cite{fs} and Strichartz\cite{s}) 
says that there exists a positive constant 
$c$ such that for any $f \in L^{2}({\mathbb R})$ and $a$ and $b$ 
real numbers, we have:
\begin{eqnarray}
\int_{{\mathbb R}}(x - a)^{2}|f(x)|^{2}dx \cdot 
\int_{{\mathbb R}}(\gamma - b)^{2}|\hat{f}(\gamma)|^{2}d\gamma 
& \geq & c\parallel f \parallel_{2}^{4}. \label{cheisenberg}
\end{eqnarray}
In the above inequality, $\hat{f}$ denotes the Fourier transform 
of $f$:
\[
\hat{f}(\gamma)=\int_{\mathbf{R}}e^{-2\pi i\gamma x}f(x)dx
\]
while $\parallel f \parallel_{2}$ represents the $L^{2}$-norm of $f$. 
Besides, both integrals from the left-hand side are assumed to be finite.\\
\par The mathematical interpretation of inequality (\ref{cheisenberg}) 
is that an $L^{2}$-function 
and its Fourier transform cannot both be localized.
See Folland and Sitaram\cite{fs}.
\\
\par In this paper the function $f$ will be replaced by a complex valued 
stochastic process $\{Y_{t}\}_{t \in I}$, where $I$ is an interval of the form 
$[0, T]$, for some $T > 0$, or $I=[0, \infty)$. The integration with respect 
to the Lebesgue measure will be replaced by the integration with respect to   
a real valued square integrable continuous martingale 
$\{X_{t}\}_{t \in I}$ whose 
quadratic variation $\{\langle X \rangle_{t}\}_{t \in I}$ is 
deterministic. The real numbers $a$ and $b$ will be replaced by 
two real valued deterministic (measurable) functions $g$ and $\tilde{g}$. 
The Fourier transform will be replaced by a unitary operator $\mathcal{G}$ 
that  will be defined in the next section.

\section{The ${\mathcal G}$-transform}

Let $(\Omega, {\mathcal F}, P)$ be a probability space and 
$\{{\mathcal F}_{t}\}_{t \in I}$ a filtration over ${\mathcal F}$.
Let $\{X_{t}\}_{t \in I}$, where $I = [0, T]$ or $I = [0, \infty)$, 
be a real valued square integrable continuous martingale defined on 
$\Omega$, adapted to $\{{\mathcal F}_{t}\}_{t \in I}$,  
whose quadratic variation process 
$\{\langle X \rangle_{t}\}_{t \in I}$ is deterministic. That means, 
for all $t \in I$, $\langle X \rangle_{t}$ is constant a.s.(almost surely).
Thus, as a process, $X_t$ is equal in distribution to a time-changed Brownian 
motion $B_{h(t)}$ with a deterministic time change $h(t)$.  We assume that 
$X_{0} = 0$ a.s.. Let 
$f : {\mathbb R}^{2} \to {\mathbb R}$ 
be a twice differentiable function with continuous second order partial 
derivatives. Applying It\^{o}'s formula to the semimartingales 
$X_{t}^{1} = X_{t}$ and $X_{t}^{2} = \langle X \rangle_{t}$
we obtain:
\begin{eqnarray}
f(X_{t}, \langle X \rangle_{t}) - f(X_{0}, \langle X \rangle_{0}) 
& = & 
\int_{0}^{t}\frac{\partial f}{\partial x}(X_{s}, \langle X \rangle_{s})
dX_{s} \nonumber\\ 
& + & \int_{0}^{t}\frac{\partial f}{\partial y}(X_{s}, \langle X \rangle_{s})
d\langle X \rangle_{s} \nonumber\\ 
& + & \frac{1}{2}
\int_{0}^{t}\frac{\partial^{2} f}{\partial x^{2}}
(X_{s}, \langle X \rangle_{s})d\langle X \rangle_{s}. \label{Ito}
\end{eqnarray}
From relation (\ref{Ito}), we see that, if $f$ satisfies the 
differential equation: 
\begin{eqnarray}
\frac{1}{2}\frac{\partial^{2} f}{\partial x^{2}} + 
\frac{\partial f}{\partial y} & = & 0, \label{Differential}
\end{eqnarray}
then $f(X_{t}, \langle X \rangle_{t})$ is a local martingale 
(see Durrett\cite{d}, page 70).\\
\par The function $f_{c}(x, y) = e^{cx - \frac{1}{2}c^{2}y}$ satisfies 
equation (\ref{Differential}), for all $c \in {\mathbb C}$. Thus 
$f_{c}(X_{t}, \langle X \rangle_{t})$ is a local martingale, 
for all $c \in {\mathbb C}$. Since 
$X_{t} \stackrel{\mathcal{D}}{=} B_{h(t)}$ we have $E(e^{cX_{t}}) < \infty$, 
for all $c \in {\mathbb R}$ and $t \in I$. Thus the process 
$\mathcal{E}_{c, t} := e^{cX_{t} - \frac{1}{2}c^{2}\langle X \rangle_{t}}$ 
is a martingale ($c$ is fixed).

\begin{lemma}\label{lem:product}
For all $t \geq 0$ and all $c$ and $d$ complex numbers, we have:
\begin{eqnarray}
\mathcal{E}_{c, t} \mathcal{E}_{d, t} = e^{cd \langle X \rangle_{t}}
\mathcal{E}_{c + d, t}. \label{eq:product}
\end{eqnarray}
\end{lemma}

\begin{proof}
\begin{eqnarray*}
\mathcal{E}_{c, t} \mathcal{E}_{d, t} & = & 
e^{cX_{t} - \frac{1}{2}c^{2}\langle X \rangle_{t}}
e^{dX_{t} - \frac{1}{2}d^{2}\langle X \rangle_{t}}\\
& = & e^{cd \langle X \rangle_{t}}
e^{(c + d)X_{t} - \frac{1}{2}(c + d)^{2}\langle X \rangle_{t}}\\
& = & e^{cd \langle X \rangle_{t}}\mathcal{E}_{c + d, t}.
\end{eqnarray*}
\end{proof}

For any $t \in I$, let $V_{t}$ be the vector space spanned by all 
exponential functions 
$\{\mathcal{E}_{c, t}\}_{c \in {\mathbb C}}$. Lemma \ref{lem:product} proves 
that $V_{t}$ is closed under multiplication. The constant process $1$ belongs 
to $V_{t}$ since $1 = \mathcal{E}_{0, t}$. Since $X_{t}$ is a real 
valued function, 
we have $\overline{\mathcal{E}_{c, t}} = \mathcal{E}_{\bar{c}, t} \in V_{t}$, 
for all $c \in {\mathbb C}$. 
Thus $V_{t}$ is an algebra containing the constant functions and closed under 
conjugation.

\begin{lemma}\label{lem:scalar}
For all $s$, $t \geq 0$ and all $c$ and $d$ complex numbers we have:
\begin{eqnarray}
E\left[\mathcal{E}_{c, s}\mathcal{E}_{d, t}\right] & = &
e^{cd \langle X \rangle_{s \wedge t}}, \label{eq:scalar}
\end{eqnarray}
where $s \wedge t$ denotes the minimum between $s$ and $t$.
\end{lemma} 

\begin{proof}
Let us assume that $s \leq t$. We have:
\begin{eqnarray*}
E\left[\mathcal{E}_{c, s}\mathcal{E}_{d, t}\right] & = & 
E\left[E\left[\mathcal{E}_{c, s}\mathcal{E}_{d, t} | \mathcal{F}_{s}\right]
\right]\\
& = & E\left[\mathcal{E}_{c, s}
E\left[\mathcal{E}_{d, t} | \mathcal{F}_{s}\right]\right]\\
& = & E\left[\mathcal{E}_{c, s}\mathcal{E}_{d, s}\right]\\
& = & E\left[e^{cd \langle X \rangle_{s}}
\mathcal{E}_{c + d, s}\right]\\
& = & e^{cd \langle X \rangle_{s}}
E\left[\mathcal{E}_{c + d, s}\right]\\
& = & e^{cd \langle X \rangle_{s}}
E\left[\left[\mathcal{E}_{c + d, s} | \mathcal{F}_{0}\right]\right]\\
& = & e^{cd \langle X \rangle_{s}}
E\left[\mathcal{E}_{c + d, 0}\right]\\
& = & e^{cd \langle X \rangle_{s}}.
\end{eqnarray*}
\end{proof}

\par We define the function:
$\mathcal{G} : \cup_{t \in I}V_{t} \to \cup_{t \in I}V_{t}$, 
$\mathcal{G}\mathcal{E}_{c, t} := \mathcal{E}_{-ic, t}$. We can see that 
$\mathcal{G}$ preserves the inner product. Indeed, we have:
\begin{eqnarray*}
E\left[\mathcal{G}\mathcal{E}_{c, s}
\overline{\mathcal{G}\mathcal{E}_{d, t}}\right] 
& = & \left[\mathcal{E}_{-ic, s}\overline{\mathcal{E}_{-id, t}}\right]\\
& = & \left[\mathcal{E}_{-ic, s}\mathcal{E}_{i\bar{d}, t}\right]\\
& = & e^{(-ic)(i\bar{d})X_{s \wedge t}}\\
& = & e^{c\bar{d}X_{s \wedge t}}\\
& = & E\left[\mathcal{E}_{c, s}\mathcal{E}_{\bar{d}, t}\right]\\
& = & E\left[\mathcal{E}_{c, s}\overline{\mathcal{E}_{d, t}}\right]. 
\end{eqnarray*}
Thus $\mathcal{G}$ can be uniquely extended to a unitary operator, that 
we denote also by $\mathcal{G}$, from the Hilbert space $\mathcal{H}$ 
into itself, where $\mathcal{H}$ is the closure of the vector space 
spanned by $\{\mathcal{E}_{c, t} ~|~ c \in {\mathbb C}, t \in I\}$ 
in $L^{2}(\Omega, \mathcal{F}, P)$.

\section{The Operator of Multiplication by $X_{t}$}

Let $\{X_{t}\}_{t \in I}$ be a square integrable continuous martingale 
process adapted to the filtration $\mathcal{F}_{t}$ and having a deterministic 
quadratic variation process $\{\langle X \rangle_{t}\}_{t \in I}$. We know 
that for each $t \in I$, $X_{t} \stackrel{\mathcal{D}}{=} B_{h(t)}$.

We note the following.  If $t$ is fixed, then $h(t)$ is fixed, and $X_t$ is 
normally distributed with mean 0 and variance $h(t)$.  Thus,
\begin{eqnarray}
\label{ltwo}
e^{c|X_{t}|} & \in& L^{2}(\Omega, \mathcal{F}, P).
\end{eqnarray}

Since $|X_{t}| < e^{|X_{t}|}$, we have 
\begin{eqnarray*}
|X_{t}\mathcal{E}_{c, t}| & = & |X_{t}| \cdot 
|e^{cX_{t} - \frac{c^{2}}{2}\langle X \rangle_{t}}|\\
& \leq & e^{|X_{t}|}e^{|c| \cdot |X_{t}|}
e^{-\frac{\rm{Re}(c^{2})}{2}\langle X \rangle_{t}}\\
& = & e^{-\frac{\rm{Re}(c^{2})}{2}\langle X \rangle_{t}}e^{(|c| + 1)|X_{t}|}.
\end{eqnarray*}
Since, according to (\ref{ltwo}), 
$e^{(|c| + 1)|X_{t}|} \in L^{2}(\Omega, \mathcal{F}, P)$, we conclude that 
$X_{t}\mathcal{E}_{c, t} \in L^{2}(\Omega, \mathcal{F}, P)$, for all 
$c \in {\mathbb C}$ and $t \in I$.

\begin{lemma}
For any $t \in I$ and $c \in {\mathbb C}$, we have:
\begin{eqnarray}
X_{t}\mathcal{E}_{c, t} & = & 
\lim_{r \to 0}\frac{\mathcal{E}_{r, t} - 1}{r}\mathcal{E}_{c, t}, \quad 
{\rm in} \ L^{2}-{\rm sense}. 
\label{Xasaderviative}
\end{eqnarray}
\end{lemma}

\begin{proof}
We have:
\begin{eqnarray*}
\lim_{r \to 0}\frac{\mathcal{E}_{r, t} - 1}{r} & = & 
\lim_{r \to 0}
\frac{e^{rX_{t} - \frac{r^{2}}{2}\langle X \rangle_{t}} - 1}{r}\\
& = & \frac{d}{dr}
\left(e^{rX_{t} - \frac{r^{2}}{2}\langle X \rangle_{t}}\right)\Big|_{r = 0}\\
& = & X_{t}.
\end{eqnarray*}
Thus $X_{t}\mathcal{E}_{c, t} =  
\lim_{r \to 0}\frac{\mathcal{E}_{r, t} - 1}{r}\mathcal{E}_{c, t}$ pointwise.
By Taylor's formula with Lagrange's remainder, 
for any fixed $\omega \in \Omega$, there exists 
a real number $u(\omega)$ between $0$ and 
$rX_{t}(\omega) - \frac{r^{2}}{2}\langle X \rangle_{t}$ such that
\begin{eqnarray*}
e^{rX_{t}(\omega) - \frac{r^{2}}{2}\langle X \rangle_{t}} & = & 
1 + rX_{t}(\omega) - \frac{r^{2}}{2}\langle X \rangle_{t} + 
\frac{1}{2}e^{u(\omega)}
\left[rX_{t}(\omega) - \frac{r^{2}}{2}\langle X \rangle_{t}\right]^{2}.
\end{eqnarray*}
This relation can be rewritten as:
\begin{eqnarray*}
\frac{e^{rX_{t}(\omega) - \frac{r^{2}}{2}\langle X \rangle_{t}} - 1}{r} - 
X_{t}(\omega) & = & 
- \frac{r}{2}\langle X \rangle_{t} + 
\frac{r}{2}e^{u(\omega)}
\left[X_{t}(\omega) - \frac{r}{2}\langle X \rangle_{t}\right]^{2}.
\end{eqnarray*}
From here, using the inequalities 
$e^{u(\omega)} \leq e^{|r| \cdot |X_{t}(\omega)|}$ 
and $|X_{t}(\omega)| < e^{|X_{t}(\omega)|}$, and (\ref{ltwo}), it follows 
from Lebesgue's dominated convergence theorem, that 
$\frac{e^{rX_{t}(\omega) - \frac{r^{2}}{2}\langle X \rangle_{t}} - 1}{r} - 
X_{t}(\omega) \to 0$ in $L^{2}(\Omega, \mathcal{F}, P)$, as $r \to 0$. 
\end{proof}

\section{The Differential Operator}

\begin{definition}
For any $t \in I$ and $c \in {\mathbb C}$, we define 
$D_{t}\mathcal{E}_{c, t} := c\langle X \rangle_{t}\mathcal{E}_{c, t}$.  
\end{definition}
We extend $D_{t}$ by linearity to the vector space $V_{t}$ 
spanned by the exponential functions 
$\{\mathcal{E}_{c, t}\}_{c \in {\mathbb C}}$. We must check that $D_{t}$ is 
well-defined.

\begin{lemma}
$D_{t} : V_{t} \to V_{t}$ is a well-defined operator.
\end{lemma}

\begin{proof}
To check that $D_{t}$ is a well-defined operator from $V_{t}$ into itself 
we must check that if a function $f$ can be expressed in two different ways 
as a linear combination of exponential functions, then calculating $D_{t}f$ 
using these two linear combinations we obtain the same result. 
This reduces to checking that if $\lambda_{1}$, $\lambda_{2}$, $\dots$, 
$\lambda_{N} \in {\mathbb C}$ and $c_{1}$, $c_{2}$, $\dots$, 
$c_{N} \in {\mathbb C}$, $c_{i} \neq c_{j}$, for $i \neq j$, such that:
$$\sum_{k = 1}^{N}\lambda_{k}\mathcal{E}_{c_{k}, t} = 0,$$
then 
$$\sum_{k = 1}^{N}\lambda_{k}c_{k}\langle X \rangle_{t}
\mathcal{E}_{c_{k}, t} = 0.$$
If $\sum_{k = 1}^{N}\lambda_{k}\mathcal{E}_{c_{k}, t} = 0$, then 
for any $n \in \{0$, $1$, $\dots$, $N - 1\}$ we have:
$E\left[\left(\sum_{k = 1}^{N}\lambda_{k}\mathcal{E}_{c_{k}, t}\right)
\mathcal{E}_{n, t}\right] = 0$. This means, according to formula 
(\ref{eq:scalar}), that:
\begin{eqnarray*}
\sum_{k = 1}^{N}\lambda_{k}e^{nc_{k}\langle X \rangle_{t}} & = & 0, 
\quad \forall n \in \{0, 1, \dots, N - 1\}.
\end{eqnarray*}
The above relations represent a homogenous linear system of $N$ equations 
and $N$ unknowns $\lambda_{1}$, $\lambda_{2}$, $\dots$, 
$\lambda_{N}$, having a Vandermonde determinant. Since 
$c_{i} \neq c_{j}$, for $i \neq j$, the Vandermonde determinant of the 
above system is different from zero if $\langle X \rangle_{t} \neq 0$. 
Thus if $\langle X \rangle_{t} \neq 0$, then 
$\lambda_{1} = \lambda_{2} = \dots = \lambda_{N} = 0$. Hence, either 
$\langle X \rangle_{t} = 0$ or 
$\lambda_{1} = \lambda_{2} = \dots = \lambda_{N} = 0$. Therefore,
$\sum_{k = 1}^{N}\lambda_{k}c_{k}\langle X \rangle_{t}
\mathcal{E}_{c_{k}, t} = 0$.
\end{proof}

We call $D_{t}$ {\it the differential operator}. Its initial domain is 
$V_{t}$ which is a dense subspace of the Hilbert space $\mathcal{H}_{t}$, 
where $\mathcal{H}_{t}$ is the closure of $V_{t}$ in 
$L^{2}(\Omega, \mathcal{F}, P)$. Observe that 
$\mathcal{H}_{t} \subset L^{2}(\Omega, \sigma(X_{t}), P) \subset 
L^{2}(\Omega, \mathcal{F}_{t}, P)$, where $\sigma(X_{t})$ is the smallest 
$\sigma$-field on $\Omega$ with respect to which $X_{t}$ is measurable. 
Since $X_{t} \stackrel{\mathcal{D}}{=} B_{h(t)}$, we have 
$\mathcal{H}_{t} = L^{2}(\Omega, \sigma(X_{t}), P)$ 
(see Janson\cite{j}, page 19).\\ 
\par We denote the adjoint of $D_{t}$ by $D_{t}^{*}$.

\begin{lemma}
For any $t \in I$, $V_{t}$ is contained in the domain of $D_{t}^{*}$ and 
for any $\varphi \in V_{t}$ we have:
\begin{eqnarray}
X_{t}\varphi & = & D_{t}\varphi + D_{t}^{*}\varphi. \label{x=d+d*}
\end{eqnarray} 
\end{lemma}

\begin{proof}
Since $V_{t}$ is spanned by the exponential functions we have to check 
that for all $d \in {\mathbb C}$, $\mathcal{E}_{d, t}$ belongs to the 
domain of $D_{t}^{*}$ and for any $c \in {\mathbb C}$ we have:
\begin{eqnarray}
\langle D_{t}\mathcal{E}_{c, t}, \mathcal{E}_{d, t}\rangle & = & 
\langle \mathcal{E}_{c, t}, (X_{t} - D_{t})\mathcal{E}_{d, t}\rangle, 
\label{adjoint}
\end{eqnarray} 
where $\langle , \rangle$ denotes the $L^{2}$-inner product.
Let us observe that, for any complex number $c$,
the function $g_{c} : {\mathbb R}^{2} \to {\mathbb R}$, 
$g_{c}(x, y) = (x - cy)e^{cx - \frac{c^{2}}{2}y}$ satisfies the differential
equation:
$$\frac{1}{2}\frac{\partial^{2}g_{c}}{\partial x^{2}} + 
\frac{\partial g_{c}}{\partial y} = 0.$$ 
It follows from It\^{o}'s formula that the process 
$\{g_{c}(X_{t}, \langle X \rangle_{t})\}_{t \in I}$ is a martingale. 
We denote this martingale process by 
$\{\mathcal{X}_{t}\mathcal{E}_{c, t}\}_{t \in I}$. So, 
$\mathcal{X}_{t}\mathcal{E}_{c, t} = (X_{t} - c\langle X \rangle_{t})
e^{cX_{t} - \frac{c^{2}}{2}\langle X \rangle_{t}}$. From this definition 
we obtain the formula:
\begin{eqnarray}
X_{t}\mathcal{E}_{c, t} & = & \mathcal{X}_{t}\mathcal{E}_{c, t} + 
c\langle X \rangle_{t}\mathcal{E}_{c, t}. \label{madjoint}
\end{eqnarray}  
Thus we have:
\begin{eqnarray*} 
\langle \mathcal{E}_{c, t}, (X_{t} - D_{t})\mathcal{E}_{d, t}\rangle 
& = & E[X_{t}\mathcal{E}_{c, t}\mathcal{E}_{\bar{d}, t}] - 
E[\mathcal{E}_{c, t}D_{t}\mathcal{E}_{\bar{d}, t}]\\
& = & E[X_{t}\mathcal{E}_{c + \bar{d}, t}e^{c\bar{d}\langle X \rangle_{t}}] 
- E[\mathcal{E}_{c, t}\bar{d}\langle X \rangle_{t}\mathcal{E}_{\bar{d}, t}]
\\
& = & e^{c\bar{d}\langle X \rangle_{t}}E[X_{t}\mathcal{E}_{c + \bar{d}, t}]
- \bar{d}\langle X \rangle_{t}E[\mathcal{E}_{c, t}\mathcal{E}_{\bar{d}, t}].
\end{eqnarray*}
Using relations (\ref{madjoint}) and (\ref{eq:scalar}) we obtain:
\begin{eqnarray*} 
& \ & \langle \mathcal{E}_{c, t}, (X_{t} - D_{t})\mathcal{E}_{d, t}\rangle\\
& = & e^{c\bar{d}\langle X \rangle_{t}}
E[\mathcal{X}_{t}\mathcal{E}_{c + \bar{d}, t} + 
(c + \bar{d})\langle X \rangle_{t}\mathcal{E}_{c + \bar{d}, t}] -  
\bar{d}\langle X \rangle_{t}E[\mathcal{E}_{c, t}\mathcal{E}_{\bar{d}, t}]\\
& = & e^{c\bar{d}\langle X \rangle_{t}}
E[\mathcal{X}_{t}\mathcal{E}_{c + \bar{d}, t}] + 
(c + \bar{d})\langle X \rangle_{t}e^{c\bar{d}\langle X \rangle_{t}}
E[\mathcal{E}_{c + \bar{d}, t}] - 
\bar{d}\langle X \rangle_{t}e^{c\bar{d}\langle X \rangle_{t}}.
\end{eqnarray*}
Because $\{\mathcal{X}_{t}\mathcal{E}_{c + \bar{d}, t}\}_{t \in I}$ and 
$\{\mathcal{E}_{c + \bar{d}, t}\}_{t \in I}$ are martingale processes 
we have: 
\begin{eqnarray*} 
E[\mathcal{X}_{t}\mathcal{E}_{c + \bar{d}, t}] & = & 
E[\mathcal{X}_{0}\mathcal{E}_{c + \bar{d}, 0}]\\
& = & 0
\end{eqnarray*} 
and 
\begin{eqnarray*}
E[\mathcal{E}_{c + \bar{d}, t}] & = & E[\mathcal{E}_{c + \bar{d}, 0}]\\
& = & 1.
\end{eqnarray*} 
Thus, we obtain:
\begin{eqnarray*}
& \ & \langle \mathcal{E}_{c, t}, (X_{t} - D_{t})\mathcal{E}_{d, t}\rangle\\
& = & e^{c\bar{d}\langle X \rangle_{t}} \cdot 0 + 
(c + \bar{d})\langle X \rangle_{t}e^{c\bar{d}\langle X \rangle_{t}} \cdot 1
- \bar{d}\langle X \rangle_{t}e^{c\bar{d}\langle X \rangle_{t}}\\
& = & c\langle X \rangle_{t}e^{c\bar{d}\langle X \rangle_{t}}\\
& = & c\langle X \rangle_{t}E[\mathcal{E}_{c, t}\mathcal{E}_{\bar{d}, t}]\\
& = & \langle c\langle X \rangle_{t}\mathcal{E}_{c, t}, 
\mathcal{E}_{d, t}\rangle\\
& = & \langle D_{t}\mathcal{E}_{c, t}, \mathcal{E}_{d, t}\rangle.
\end{eqnarray*}
\end{proof}

\begin{corollary}
For all $t \in I$, the differential operator $D_{t}$ 
admits a closed extension. 
\end{corollary}

\begin{proof} Let $W_{t}$ be the subspace of all functions $\varphi$
in $\mathcal{H}_{t}$ for which there exists a sequence
$\{\varphi_{n}\}_{n \geq 1} \subset V_{t}$ such that     
$\varphi_{n} \to \varphi$ in $L^{2}(\Omega, \mathcal{F}, P)$ and the 
sequence $\{D_{t}\varphi_{n}\}_{n \geq 1}$ is Cauchy in
$L^{2}(\Omega, \mathcal{F}, P)$. 
We define the operator $\widetilde{D}_{t}$ on $W_{t}$ in the following way: 
if $\varphi \in W_{t}$ and $\{\varphi_{n}\}_{n \geq 1} \subset V_{t}$ such 
that $\varphi_{n} \to \varphi$ in $ L^{2}(\Omega, \mathcal{F}, P)$ and 
$\{D_{t}\varphi_{n}\}_{n \geq 1}$ is Cauchy in 
$L^{2}(\Omega, \mathcal{F}, P)$, then $\widetilde{D}_{t}\varphi = u$, where 
$u = \lim_{n \to \infty}D_{t}\varphi_{n}$ in $L^{2}(\Omega, \mathcal{F}, P)$. 
We notice that $u \in \mathcal{H}_{t}$.
\par We need to check that $\widetilde{D}_{t}\varphi$ is well defined. 
Let $\{\varphi_{n}\}_{n \geq 1} \subset V_{t}$ such 
that $\varphi_{n} \to \varphi$ in $ L^{2}(\Omega, \mathcal{F}, P)$ and 
$\{D_{t}\varphi_{n}\}_{n \geq 1}$ is Cauchy in 
$L^{2}(\Omega, \mathcal{F}, P)$, and 
$\{\psi_{n}\}_{n \geq 1} \subset V_{t}$ such 
that $\psi_{n} \to \varphi$ in $ L^{2}(\Omega, \mathcal{F}, P)$ and 
$\{D_{t}\psi_{n}\}_{n \geq 1}$ is Cauchy in 
$L^{2}(\Omega, \mathcal{F}, P)$. 
Let $u := \lim_{n \to \infty}D_{t}\varphi_{n}$ and 
$v := \lim_{n \to \infty}D_{t}\psi_{n}$. Both limits are in the $L^{2}$-sense. 
For any $\phi \in V_{t}$ we have:
\begin{eqnarray*}
\langle u, \phi \rangle & = & 
\lim_{n \to \infty}\langle D_{t}\varphi_{n}, \phi \rangle\\
& = & \lim_{n \to \infty}\langle \varphi_{n}, D_{t}^{*}\phi \rangle\\
& = & \langle \varphi, D_{t}^{*}\phi \rangle.
\end{eqnarray*}
In the same way, we can see that
\begin{eqnarray*}
\langle v, \phi \rangle & = & \langle \varphi, D_{t}^{*}\phi \rangle.
\end{eqnarray*}
Thus $\langle u, \phi \rangle = \langle v, \phi \rangle$, for all 
$\phi \in V_{t}$. Since $V_{t}$ is dense in $\mathcal{H}_{t}$ and 
$u$, $v \in \mathcal{H}_{t}$ it follows that $u = v$.
\end{proof}

From now on we will denote the closure of $D_{t}$ by $D_{t}$, too.

\section{Commutation Relationships}

\begin{lemma}
For any $t \in I$, we have: 
\begin{eqnarray}
[D_{t}, X_{t}] & = & \langle X \rangle_{t}I. 
\end{eqnarray}
\end{lemma}

\begin{proof}
We will check that for any exponential function $\mathcal{E}_{c, t}$ we have:
\begin{eqnarray}
D_{t}X_{t}\mathcal{E}_{c, t} - X_{t}D_{t}\mathcal{E}_{c, t} & = & 
\langle X \rangle_{t}\mathcal{E}_{c, t}. \label{comdx}
\end{eqnarray}
Indeed we have:
\begin{eqnarray*}
D_{t}X_{t}\mathcal{E}_{c, t} & = & D_{t}\lim_{s \to 0}
\frac{\mathcal{E}_{s, t} - 1}{s}\mathcal{E}_{c, t}.
\end{eqnarray*}
We have seen that, because of (\ref{ltwo}), the limit from 
the right-hand side of the last equality is both pointwise and in the 
$L^{2}$-sense. The last equality can be written now as:
\begin{eqnarray*}
D_{t}X_{t}\mathcal{E}_{c, t} & = & D_{t}\lim_{s \to 0}
\frac{\mathcal{E}_{s, t}\mathcal{E}_{c, t} - \mathcal{E}_{c, t}}{s}\\
& = &  D_{t}\lim_{s \to 0}
\frac{\mathcal{E}_{s + c, t}
e^{sc\langle X \rangle_{t}} - \mathcal{E}_{c, t}}{s}.
\end{eqnarray*}
Because we are working with the closure of the operator $D_{t}$ and the 
limit from the last relation is in the $L^{2}$-sense, we can commute 
$D_{t}$ with the limit if $D_{t}\frac{\mathcal{E}_{s + c, t}
e^{sc\langle X \rangle_{t}} - \mathcal{E}_{c, t}}{s}$ converges in 
$L^{2}(\Omega, \mathcal{F}, P)$, as $s \to 0$. 
Commuting $D_{t}$ with the limit we get:
\begin{eqnarray*}
D_{t}X_{t}\mathcal{E}_{c, t}
& = & \lim_{s \to 0}D_{t}\frac{\mathcal{E}_{s + c, t}
e^{sc\langle X \rangle_{t}} - \mathcal{E}_{c, t}}{s}\\
& = & \lim_{s \to 0}\frac{e^{sc\langle X \rangle_{t}}
D_{t}\mathcal{E}_{s + c, t} - D_{t}\mathcal{E}_{c, t}}{s}\\
& = & \lim_{s \to 0}\frac{e^{sc\langle X \rangle_{t}}(s + c)
\langle X \rangle_{t}\mathcal{E}_{s + c, t} - 
c\langle X \rangle_{t}\mathcal{E}_{c, t}}{s}\\
& = & \langle X \rangle_{t}\lim_{s \to 0}
\frac{(s + c)e^{sc\langle X \rangle_{t}}\mathcal{E}_{s + c, t} - 
c\mathcal{E}_{c, t}}{s}\\ 
& = & \langle X \rangle_{t}\lim_{s \to 0}
\frac{(s + c)\mathcal{E}_{c, t}\mathcal{E}_{s, t} - 
c\mathcal{E}_{c, t}}{s}\\
& = &  \langle X \rangle_{t}
\left[\lim_{s \to 0}\mathcal{E}_{s, t} + c\lim_{s \to 0}
\frac{\mathcal{E}_{s, t} - 1}{s}\right]\mathcal{E}_{c, t}\\
& = & \langle X \rangle_{t}[\mathcal{E}_{0, t} + cX_{t}]
\mathcal{E}_{c, t}\\
& = & \langle X \rangle_{t}(1 + cX_{t})\mathcal{E}_{c, t}\\
& = & \langle X \rangle_{t}\mathcal{E}_{c, t} + 
X_{t}(c\langle X \rangle_{t}\mathcal{E}_{c, t})\\
& = & \langle X \rangle_{t}\mathcal{E}_{c, t} + 
X_{t}D_{t}\mathcal{E}_{c, t}.
\end{eqnarray*}
Because of (\ref{ltwo}) all the above limits are not only 
pointwise, but also in the $L^{2}$-sense.
Thus 
\begin{eqnarray*}
D_{t}X_{t}\mathcal{E}_{c, t} - X_{t}D_{t}\mathcal{E}_{c, t} = 
\langle X \rangle_{t}\mathcal{E}_{c, t}.
\end{eqnarray*}
\end{proof}

\begin{corollary}
For any $t \in I$, we have:
\begin{eqnarray}
[D_{t}, D_{t}^{*}] & = & \langle X \rangle_{t}I. \label{comdd*}
\end{eqnarray}
\end{corollary}

\begin{proof}
According to formula (\ref{x=d+d*}) we have $X_{t} = D_{t} + D_{t}^{*}$. 
Thus:
\begin{eqnarray*}
[D_{t}, D_{t}^{*}] & = & [D_{t}, X_{t} - D_{t}]\\
& = & [D_{t}, X_{t}] - [D_{t}, D_{t}]\\
& = & \langle X \rangle_{t}I - 0\\
& = & \langle X \rangle_{t}I.
\end{eqnarray*}
\end{proof}

\begin{lemma}
For any $t \in I$, we have:
\begin{eqnarray}
D_{t}\mathcal{G} & = & -i\mathcal{G}D_{t}. \label{dg} 
\end{eqnarray}
\end{lemma}

\begin{proof}
For any exponential function $\mathcal{E}_{c, t}$, we have:
\begin{eqnarray*}
D_{t}\mathcal{G}\mathcal{E}_{c, t} & = & D_{t}\mathcal{E}_{-ic, t}\\
& = & -ic\langle X \rangle_{t}\mathcal{E}_{-ic, t}\\
& = & -ic\langle X \rangle_{t}\mathcal{G}\mathcal{E}_{c, t}\\
& = & -i\mathcal{G}\left(c\langle X \rangle_{t}\mathcal{E}_{c, t}\right)\\
& = & -i\mathcal{G}D_{t}\mathcal{E}_{c, t}.
\end{eqnarray*}
\end{proof}

\begin{lemma}
For any $t \in I$, we have:
\begin{eqnarray}
D_{t}^{*}\mathcal{G} & = & i\mathcal{G}D_{t}^{*}. \label{d*g} 
\end{eqnarray}
\end{lemma}

\begin{proof}
For any two exponential functions $\mathcal{E}_{c, t}$ and 
$\mathcal{E}_{d, t}$, we have:
\begin{eqnarray*}
\langle D_{t}^{*}\mathcal{G}\mathcal{E}_{c, t}, \mathcal{E}_{d, t} \rangle
& = & 
\langle \mathcal{G}\mathcal{E}_{c, t}, D_{t}\mathcal{E}_{d, t} \rangle\\
& = & \langle \mathcal{E}_{-ic, t}, 
d\langle X \rangle_{t}\mathcal{E}_{d, t} \rangle\\
& = & \bar{d}\langle X \rangle_{t}\langle \mathcal{E}_{-ic, t}, 
\mathcal{E}_{d, t} \rangle\\ 
& = & \bar{d}\langle X \rangle_{t}e^{-ic\bar{d}\langle X \rangle_{t}}\\ 
& = & i(-i\bar{d})\langle X \rangle_{t}
e^{c\overline{(id)}\langle X \rangle_{t}}\\
& = & i(-i\bar{d})\langle X \rangle_{t}
\langle \mathcal{E}_{c, t}, \mathcal{E}_{id, t} \rangle\\
& = & i\langle \mathcal{E}_{c, t}, 
id\langle X \rangle_{t}\mathcal{E}_{id, t} \rangle\\
& = & i\langle \mathcal{E}_{c, t}, D_{t}\mathcal{E}_{id, t} \rangle\\
& = & i\langle D_{t}^{*}\mathcal{E}_{c, t}, \mathcal{E}_{id, t} \rangle.
\end{eqnarray*}
Since $\mathcal{G}$ preserves the inner product, being a unitary operator, 
we obtain:
\begin{eqnarray*}
\langle D_{t}^{*}\mathcal{G}\mathcal{E}_{c, t}, \mathcal{E}_{d, t} \rangle
& = & i\langle D_{t}^{*}\mathcal{E}_{c, t}, \mathcal{E}_{id, t} \rangle\\
& = & i\langle \mathcal{G}D_{t}^{*}\mathcal{E}_{c, t}, 
\mathcal{G}\mathcal{E}_{id, t} \rangle\\
& = & i\langle \mathcal{G}D_{t}^{*}\mathcal{E}_{c, t}, 
\mathcal{E}_{-i(id), t} \rangle\\
& = & \langle i\mathcal{G}D_{t}^{*}\mathcal{E}_{c, t}, 
\mathcal{E}_{d, t} \rangle.
\end{eqnarray*}
Thus $D_{t}^{*}\mathcal{G} = i\mathcal{G}D_{t}^{*}$.
\end{proof}

\section{Heisenberg Inequality}

\begin{lemma}
Let $t \in I$ and $Y_{t} \in L^{2}(\Omega, \sigma(X_{t}), P)$. 
Let $c_{t}$, $\tilde{c}_{t} \in {\mathbb R}$. 
If $Y_{t}$ belongs to the domains of the operators $D_{t}$, $D_{t}^{*}$, 
$D_{t}D_{t}^{*}$, and $D_{t}^{*}D_{t}$, then 
\begin{eqnarray}
E[(X_{t} - c_{t})^{2}|Y_{t}|^{2}]^{\frac{1}{2}}
E[(X_{t} - \tilde{c}_{t})^{2}|\mathcal{G}Y_{t}|^{2}]^{\frac{1}{2}} & \geq & 
\langle X \rangle_{t}E[|Y_{t}|^{2}]. \label{h1}  
\end{eqnarray}
\end{lemma}

\begin{proof}
Since $X_{t} = D_{t} + D_{t}^{*}$, we have:
\begin{eqnarray*}
& \ & E[(X_{t} - c_{t})^{2}|Y_{t}|^{2}]^{\frac{1}{2}}
E[(X_{t} - \tilde{c}_{t})^{2}|\mathcal{G}Y_{t}|^{2}]^{\frac{1}{2}}\\
& = & E[|(D_{t} + D_{t}^{*} - c_{t})Y_{t}|^{2}]^{\frac{1}{2}}
E[|(D_{t} + D_{t}^{*} - \tilde{c}_{t})\mathcal{G}Y_{t}|^{2}]^{\frac{1}{2}}\\
& = & E[|(D_{t} + D_{t}^{*} - c_{t})Y_{t}|^{2}]^{\frac{1}{2}}
E[|(D_{t}\mathcal{G} + D_{t}^{*}\mathcal{G} - \tilde{c}_{t}\mathcal{G})
Y_{t}|^{2}]^{\frac{1}{2}}.
\end{eqnarray*}
Using the commutation relationships (\ref{dg}) and (\ref{d*g}) we obtain:
\begin{eqnarray*}
& \ & E[(X_{t} - c_{t})^{2}|Y_{t}|^{2}]^{\frac{1}{2}}
E[(X_{t} - \tilde{c}_{t})^{2}|\mathcal{G}Y_{t}|^{2}]^{\frac{1}{2}}\\
& = & E[|(D_{t} + D_{t}^{*} - c_{t})Y_{t}|^{2}]^{\frac{1}{2}}
E[|(-i\mathcal{G}D_{t} + i\mathcal{G}D_{t}^{*} - \tilde{c}_{t}\mathcal{G})
Y_{t}|^{2}]^{\frac{1}{2}}\\
& = & E[|(D_{t} + D_{t}^{*} - c_{t})Y_{t}|^{2}]^{\frac{1}{2}}
E[|i\mathcal{G}(-D_{t} + D_{t}^{*} + i\tilde{c}_{t})Y_{t}|^{2}]^{\frac{1}{2}}.
\end{eqnarray*}
Since $\mathcal{G}$ is a unitary operator and $i$ has modulus 1, we have 
$E[|i\mathcal{G}(-D_{t} + D_{t}^{*} + i\tilde{c}_{t})Y_{t}|^{2}] = 
E[|(-D_{t} + D_{t}^{*} + i\tilde{c}_{t})Y_{t}|^{2}]$. 
Thus we obtain:
\begin{eqnarray*}
& \ & E[(X_{t} - c_{t})^{2}|Y_{t}|^{2}]^{\frac{1}{2}}
E[(X_{t} - \tilde{c}_{t})^{2}|\mathcal{G}Y_{t}|^{2}]^{\frac{1}{2}}\\
& = & E[|(D_{t} + D_{t}^{*} - c_{t})Y_{t}|^{2}]^{\frac{1}{2}}
E[|(-D_{t} + D_{t}^{*} + i\tilde{c}_{t})Y_{t}|^{2}]^{\frac{1}{2}}\\
& = & \frac{1}{2}\left\{E[|(D_{t} + D_{t}^{*} - c_{t})Y_{t}|^{2}]^{\frac{1}{2}}
E[|(-D_{t} + D_{t}^{*} + i\tilde{c}_{t})Y_{t}|^{2}]^{\frac{1}{2}}\right.\\
& + & 
\left.E[|(-D_{t} + D_{t}^{*} + i\tilde{c}_{t})Y_{t}|^{2}]^{\frac{1}{2}}
E[|(D_{t} + D_{t}^{*} - c_{t})Y_{t}|^{2}]^{\frac{1}{2}}\right\}.
\end{eqnarray*}
Applying Schwarz inequality we obtain:
\begin{eqnarray*}
& \ & E[(X_{t} - c_{t})^{2}|Y_{t}|^{2}]^{\frac{1}{2}}
E[(X_{t} - \tilde{c}_{t})^{2}|\mathcal{G}Y_{t}|^{2}]^{\frac{1}{2}}\\
& = & \frac{1}{2}\left\{E[|(D_{t} + D_{t}^{*} - c_{t})Y_{t}|^{2}]^{\frac{1}{2}}
E[|(-D_{t} + D_{t}^{*} + i\tilde{c}_{t})Y_{t}|^{2}]^{\frac{1}{2}}\right.\\
& + & 
\left.E[|(-D_{t} + D_{t}^{*} + i\tilde{c}_{t})Y_{t}|^{2}]^{\frac{1}{2}}
E[|(D_{t} + D_{t}^{*} - c_{t})Y_{t}|^{2}]^{\frac{1}{2}}\right\}.\\
& \geq & \frac{1}{2}\left\{\left|\langle(D_{t} + D_{t}^{*} - c_{t})Y_{t}, 
(-D_{t} + D_{t}^{*} + i\tilde{c}_{t})Y_{t}\rangle\right|\right.\\
& + & \left.\left|\langle(-D_{t} + D_{t}^{*} + i\tilde{c}_{t})Y_{t}, 
(D_{t} + D_{t}^{*} - c_{t})Y_{t}\rangle\right|\right\}.
\end{eqnarray*}
Applying the triangle inequality we get:
\begin{eqnarray*}
& \ & E[(X_{t} - c_{t})^{2}|Y_{t}|^{2}]^{\frac{1}{2}}
E[(X_{t} - \tilde{c}_{t})^{2}|\mathcal{G}Y_{t}|^{2}]^{\frac{1}{2}}\\
& \geq & \frac{1}{2}\left|\langle(D_{t} + D_{t}^{*} - c_{t})Y_{t}, 
(-D_{t} + D_{t}^{*} + i\tilde{c}_{t})Y_{t}\rangle\right.\\
& + & \left.\langle(-D_{t} + D_{t}^{*} + i\tilde{c}_{t})Y_{t}, 
(D_{t} + D_{t}^{*} - c_{t})Y_{t}\rangle\right|.
\end{eqnarray*}
Using the linearity in the first argument and  
the conjugate linearity in the second argument of the inner product, and the 
definition of the adjoint operator, because $c_{t}$ and $\tilde{c}_{t}$ 
are real numbers, the last expression can be simplified and we obtain 
finally:
\begin{eqnarray*}
& \ & E[(X_{t} - c_{t})^{2}|Y_{t}|^{2}]^{\frac{1}{2}}
E[(X_{t} - \tilde{c}_{t})^{2}|\mathcal{G}Y_{t}|^{2}]^{\frac{1}{2}}\\
& \geq & 
\frac{1}{2}\left|-2\langle D_{t}Y_{t}, D_{t}Y_{t}\rangle
 + 2\langle D_{t}^{*}Y_{t}, D_{t}^{*}Y_{t}\rangle\right|\\
& = & 
\left|-\langle D_{t}Y_{t}, D_{t}Y_{t}\rangle
 + \langle D_{t}^{*}Y_{t}, D_{t}^{*}Y_{t}\rangle\right|\\
& = & 
\left|-\langle D_{t}^{*}D_{t}Y_{t}, Y_{t}\rangle
 + \langle D_{t}D_{t}^{*}Y_{t}, Y_{t}\rangle\right|\\
& = & 
\left|\langle (D_{t}D_{t}^{*} - D_{t}^{*}D_{t})Y_{t}, Y_{t}\rangle\right|\\
& = & 
\left|\langle [D_{t}, D_{t}^{*}]Y_{t}, Y_{t}\rangle\right|.
\end{eqnarray*}
Using the commutation relationship (\ref{comdd*}) we obtain:
\begin{eqnarray*}
& \ & E[(X_{t} - c_{t})^{2}|Y_{t}|^{2}]^{\frac{1}{2}}
E[(X_{t} - \tilde{c}_{t})^{2}|\mathcal{G}Y_{t}|^{2}]^{\frac{1}{2}}\\
& \geq & 
\left|\langle \langle X \rangle_{t}Y_{t}, Y_{t}\rangle\right|\\
& = & \langle X \rangle_{t}\langle Y_{t}, Y_{t}\rangle\\
& = & \langle X \rangle_{t}E[|Y_{t}|^{2}].
\end{eqnarray*}
\end{proof}

\begin{theorem}
Let $\{Y_{t}\}_{t \in [0, T]}$ be a stochastic process 
such that for each $t \in [0$, $T]$, 
$Y_{t} \in L^{2}(\Omega, \sigma(X_{t}), P)$ and 
$Y_{t}$ belongs to the domain of $D_{t}$, $D_{t}^{*}$, 
$D_{t}D_{t}^{*}$, and $D_{t}^{*}D_{t}$. 
Then for any Borel measurable functions
$g$, $\tilde{g} : [0, T] \to {\mathbb R}$, we have:
\begin{eqnarray}
& \ & E\left[\left|\int_{0}^{T}(X_{t} - g(t))Y_{t}dX_{t}
\right|^{2}\right]^{\frac{1}{2}}
E\left[\left|\int_{0}^{T}(X_{t} - \tilde{g}(t))
\mathcal{G}Y_{t}dX_{t}\right|^{2}\right]^{\frac{1}{2}} \nonumber\\
& \geq & E\left[\int_{0}^{T}|Y_{t}|^{2}\langle X \rangle_{t}
d\langle X \rangle_{t}\right]. \label{h2}
\end{eqnarray}
We assume that all the measurability and square integrability conditions, 
necessary for the existence of the above stochastic integrals, hold.
\end{theorem}

\begin{proof} 
Using the isomorphism through which the stochastic integral is extended 
from the simple processes to the square integrable processes we obtain:
\begin{eqnarray*}
& \ & E\left[\left|\int_{0}^{T}(X_{t} - g(t))Y_{t}dX_{t}
\right|^{2}\right]
E\left[\left|\int_{0}^{T}(X_{t} - \tilde{g}(t))
\mathcal{G}Y_{t}dX_{t}\right|^{2}\right]\\
& = & 
\int_{0}^{T}E\left[|(X_{t} - g(t))Y_{t}|^{2}\right]
d\langle X \rangle_{t}
\int_{0}^{T}E\left[|(X_{t} - \tilde{g}(t))\mathcal{G}Y_{t}|^{2}\right]
d\langle X \rangle_{t}.
\end{eqnarray*}
Applying Schwarz inequality we obtain:
\begin{eqnarray*}
& \ & E\left[\left|\int_{0}^{T}(X_{t} - g(t))Y_{t}dX_{t}
\right|^{2}\right]
E\left[\left|\int_{0}^{T}(X_{t} - \tilde{g}(t))
\mathcal{G}Y_{t}dX_{t}\right|^{2}\right]\\
& = & 
\int_{0}^{T}E\left[|(X_{t} - g(t))Y_{t}|^{2}\right]
d\langle X \rangle_{t}
\int_{0}^{T}E\left[|(X_{t} - \tilde{g}(t))\mathcal{G}Y_{t}|^{2}\right]
d\langle X \rangle_{t}\\
& \geq & \left\{\int_{0}^{T}
E\left[|(X_{t} - g(t))Y_{t}|^{2}\right]^{\frac{1}{2}}
E\left[|(X_{t} - \tilde{g}(t))\mathcal{G}Y_{t}|^{2}\right]^{\frac{1}{2}}
d\langle X \rangle_{t}\right\}^{2}.
\end{eqnarray*}
Applying inequality (\ref{h1}), for $c_{t} := g(t)$ and 
$\tilde{c}_{t} := \tilde{g}(t)$, we obtain:
\begin{eqnarray*}
& \ & E\left[\left|\int_{0}^{T}(X_{t} - g(t))Y_{t}dX_{t}
\right|^{2}\right]^{\frac{1}{2}}
E\left[\left|\int_{0}^{T}(X_{t} - \tilde{g}(t))
\mathcal{G}Y_{t}dX_{t}\right|^{2}\right]^{\frac{1}{2}}\\
& \geq & \int_{0}^{T}
E\left[|(X_{t} - g(t))Y_{t}|^{2}\right]^{\frac{1}{2}}
E\left[|(X_{t} - \tilde{g}(t))\mathcal{G}Y_{t}|^{2}\right]^{\frac{1}{2}}
d\langle X \rangle_{t}\\
& \geq & \int_{0}^{T}\langle X \rangle_{t}E[|Y_{t}|^{2}]
d\langle X \rangle_{t}\\
& = & E\left[\int_{0}^{T}|Y_{t}|^{2}\langle X \rangle_{t}
d\langle X \rangle_{t}\right].
\end{eqnarray*}
\end{proof}
If we take, $Y_{t} := \varphi(X_{t})$, where 
$\varphi : {\mathbb R} \to {\mathbb C}$ is a Borel measurable function, 
then we obtain the following:

\begin{corollary}
Let $\varphi : {\mathbb R} \to {\mathbb C}$ be a Borel measurable function
such that for each $t \in [0$, $T]$,  
$\varphi(X_{t})$ belongs to the domain of $D_{t}$, $D_{t}^{*}$, 
$D_{t}D_{t}^{*}$, and $D_{t}^{*}D_{t}$. 
Then for any Borel measurable functions
$g$, $\tilde{g} : [0, T] \to {\mathbb R}$, we have:
\begin{eqnarray}
& \ & E\left[\left|\int_{0}^{T}(X_{t} - g(t))\varphi(X_{t})dX_{t}
\right|^{2}\right]^{\frac{1}{2}}
E\left[\left|\int_{0}^{T}(X_{t} - \tilde{g}(t))
\mathcal{G}\varphi(X_{t})dX_{t}\right|^{2}\right]^{\frac{1}{2}} \nonumber\\
& \geq & E\left[\int_{0}^{T}|\varphi(X_{t})|^{2}\langle X \rangle_{t}
d\langle X \rangle_{t}\right].
\end{eqnarray}
We assume that all the measurability and square integrability conditions, 
necessary for the existence of the above stochastic integrals, hold.
\end{corollary}

Suppose $X_{t} = B_{t}$ is a Brownian motion process. Then for each 
$t \in [0$, $T]$, $\langle X \rangle_{t} = t$ is deterministic. 
The inequality (\ref{h2}) becomes:

\begin{theorem}
Let $\{Y_{t}\}_{t \in [0, T]}$ be a stochastic process 
such that for each $t \in [0$, $T]$, 
$Y_{t}$ is measurable with respect to the sigma-field $\sigma(B_{t})$ and 
$Y_{t}$ belongs to the domain of $D_{t}$, $D_{t}^{*}$, 
$D_{t}D_{t}^{*}$, and $D_{t}^{*}D_{t}$. 
Then for any Borel measurable functions
$g$, $\tilde{g} : [0, T] \to {\mathbb R}$, we have:
\begin{eqnarray}
& \ & E\left[\left|\int_{0}^{T}(B_{t} - g(t))Y_{t}dB_{t}
\right|^{2}\right]^{\frac{1}{2}}
E\left[\left|\int_{0}^{T}(B_{t} - \tilde{g}(t))
\mathcal{G}Y_{t}dB_{t}\right|^{2}\right]^{\frac{1}{2}} \nonumber\\
& \geq & E\left[\int_{0}^{T}t|Y_{t}|^{2}
dt\right].
\end{eqnarray}
We assume that all the measurability and square integrability conditions, 
necessary for the existence of the above stochastic integrals, hold.
\end{theorem}

In all the inequalities from this section we may replace $\int_{0}^{T}$ by 
$\int_{0}^{\infty}$ and obtain the same results.

\bibliographystyle{plain}

\end{document}